\documentclass[11pt]{amsart}
\usepackage{amsmath}
\usepackage{amssymb}
\usepackage{amsthm}
\usepackage{latexsym}

\newtheorem{Theorem}{Theorem}[section]
\newtheorem{Lemma}[Theorem]{Lemma}
\newtheorem{Corollary}[Theorem]{Corollary}
\newtheorem{Proposition}[Theorem]{Proposition}
\newtheorem{Remark}[Theorem]{Remark}

\newtheorem{Definition}[Theorem]{Definition}

\newtheorem{Notation}[Theorem]{Notation}

\def\opn#1#2{\def#1{\operatorname{#2}}}
\opn\depth{depth}
\opn\reg{reg}
\opn\height{height}
\opn\pd{pd}
\opn\Tor{Tor}
\opn\Ass{Ass}
\opn\multideg{multideg}

\numberwithin{equation}{section}

\begin{document}

\title{Some algebraic invariants related to mixed product ideals}
\author{Cristodor IONESCU  and Giancarlo RINALDO }
\address{Cristodor Ionescu, Institute of Mathematics Simion Stoilow of the Romanian Academy, P.O. Box 1-764, 014700 Bucharest. Romania}
\email{cristodor.ionescu@imar.ro}
\thanks{The first author was  supported by the  CNCSIS grant ID-PCE no. 51/2007}
\address{Giancarlo Rinaldo, Dipartimento di Matematica, Universit\`a di Messina , Contrada
Papardo, salita Sperone, 31, 98166 Messina. Italy}
\email{rinaldo@dipmat.unime.it}

\begin{abstract}

We compute some algebraic invariants (e.g. depth, Castelnuovo - Mumford regularity) for a special class of monomial ideals, namely the ideals of mixed products. As a consequence, we characterize the Cohen-Macaulay ideals of mixed products.
\end{abstract}

\maketitle

\section{Introduction and preliminaries}
\label{1}

 The class of ideals of mixed products is a special class of square-free monomial ideals. They were first introduced by Restuccia and Villarreal (see \cite{RV} and \cite{Vi}), who studied the normality of such ideals. They  gave a complete classification of normal mixed products ideals, as well as applications in graph theory. 

Let $S=K[\mathbf{x},\mathbf{y}]$ be a polynomial ring over a field $K$ in two disjoint
sets of variables $\mathbf{x}=(x_1,\dots,x_n)$, $\mathbf{y}=(y_1,\dots,y_m)$. The \textit{ideals of mixed products} are the ideals 
$$I_qJ_r+I_sJ_t,\ \ q,r,s,t\in{\bf N},\ \ q+r=s+t,$$
where $I_l$ (resp. $J_p$) is the ideal of $S$ generated by all the
square-free monomials of degree $l$ (resp. $p$) in the variables $\mathbf{x}$
(resp. $\mathbf{y}$). We set  $I_0=J_0=S$. 
By symmetry, essentially there are  2 cases:

   $$ i)\ L=I_kJ_r+I_sJ_t,\ 0\leq k< s,$$
  $$  ii)\ L=I_kJ_r,\ k\geq 1\ {\rm or}\ r\geq 1.$$ 
Our aim is to investigate some algebraic invariants of this type of ideals, such as depth and Castelnuovo-Mumford regularity. One case is already known. Namely in \cite{HH},  Herzog and Hibi proved the case $m=0$ (see Proposition \ref{lem_reg}). In the general case, we use also some techniques used by Herzog, Restuccia and the second author in the paper \cite{HRR}. 
Note that our results are true in a slightly more general situation than the one considered in \cite{RV}, namely we don't use the condition $q+r=s+t.$

\par\noindent All over the paper we shall use the following notation. 

\begin{Notation}\label{notatie} {\rm Let $K$ be a field and $K[\mathbf{x},\mathbf{y}]=K[x_1,\dots,x_n,y_1,\dots,y_m]$ the polynomial ring in $n+m$ variables over $K.$ By $I_k$ we shall mean the  monomial ideal generated by all the square-free monomials of degree $k$ in the variables  $x_1,\dots,x_n$ and by $J_l$ we shall mean the  monomial ideal generated by all the square-free monomials of degree $l$ in the variables $y_1,\dots,y_m$.
 For a monomial ideal $I$, we shall denote by $G(I)$ the minimal monomial system of generators of $I.$}
\end{Notation}
We recall the following
\begin{Definition}
  Let $I=(\mathbf{x}^{\alpha_1},\ldots, \mathbf{x}^{\alpha_q})\in K[\mathbf{x}]= K[x_1,\ldots,x_n]$ be a square-free
  monomial ideal, with $\alpha_i=(\alpha_{i_1},\ldots,\alpha_{i_n})\in \{0,1\}^n$. The \emph{Alexander dual}
  of $I$ is the ideal
  
  $$I^*=\mathop\bigcap\limits_{i=1}^q \mathfrak{m}_i,$$
  
  where $\mathfrak{m}_i=(x_j:\alpha_{i_j}=1)$.
\end{Definition}

\begin{Remark}\label{rem:mixdual}
  {\rm It is easy to see that 
  $$I_k=I^*_{n-k+1},\ k=1,\ldots,n.$$ 
 It follows that 
  $$I_k=\cap (x_{i_1},\ldots, x_{i_{n-k+1}}),$$
  where the intersection is taken over all subsets with $n-k+1$ elements $\{i_1,\dots ,i_{n-k+1}\}\subseteq\{1,\dots,n\}$. Hence
  $$\ \dim (S/I_k)=k-1.$$}
\end{Remark}

\section{Regularity of ideals of mixed products}
\label{sec:regularity}
In this section we want to study the Castelnuovo-Mumford regularity of  
the mixed product ideal  $I_q J_r+I_s J_t$.

We remind that, if $\mathbb{F}$ is the minimal graded free resolution of a given graded finitely generated $S-$module $M$
$$\mathbb{F}:\  0\to \mathop\oplus\limits_{i=1}^{b_g}S(-d_{g_i}) \stackrel{\phi_g}{\to} \cdots \to
\mathop\oplus\limits_{i=1}^{b_k}S(-d_{k_i})
\stackrel{\phi_k}{\to}\cdots\to\mathop\oplus\limits_{i=1}^{b_0} S(-d_{0_i})\stackrel{\phi_0}{\to} M\rightarrow 0,$$
then the \textit{Castelnuovo-Mumford regularity of $M$} is given by 
$$\reg(M):=\max (d_{k_i}-k\ |\ k >0, i=1,\ldots,b_k ).$$
If  $J\subseteq K[z_1,\ldots,z_s]$ is a graded ideal, if all the generators of $J$ have the same degree $d$ and also $\reg J=d$, then 
$J$ has a \textit{d-linear free resolution}.

We recall the following (see for example \cite[Th. 5.56]{MS}).
\begin{Theorem}[Eagon-Reiner Theorem]
  Let $I$ be a square-free monomial ideal. Then $S/I$ is
  Cohen-Macaulay if and only if $I^*$ has a linear free resolution.
\end{Theorem}

\begin{Proposition}\label{lem_reg}
Suppose that $m=0$, so that $S=K[x_1,\ldots,x_n]$. Then:
\par\noindent {\rm a)} $ S/I_k$ is Cohen-Macaulay for every $k$;
\par\noindent {\rm b)} $\reg (I_k)=k$.
\par\noindent {\rm c)} $\pd (S/I_k)=n-k+1.$
\end{Proposition}

\textit{Proof}: a) See \cite[Ex. 2.2]{HH}.
\par\noindent b)  By Remark \ref{rem:mixdual} we have $\dim S/I_k=k-1$. Therefore, by \cite[Theorem 5.59]{MS} and a) we get
  $$\reg(I_k)=\pd(S_n/I^*_k)=\pd(S_n/I_{n-k+1})=k.$$  
\par\noindent c) $\pd (S/I_k)=\reg(S/I^*_k)=\reg(S/I_{n-k+1})= n-k+1.$

\begin{Lemma}
  Let $I$ and $J$ be graded ideals in $S$ such that
  $\Tor_1^S(S/I,S/J)=0.$ Assume that $\reg(I)=q$  and $\reg(J)=r.$
   Then \[\reg(I+J)=q+r-1.\]
\end{Lemma}

\textit{Proof}:  Let $\mathbb{F}=(F_i)$ and $\mathbb{G}=(G_j)$ be the graded minimal free resolutions of $S/I,$ resp. $S/J$  and consider the tensor product of  complexes $\mathbb{F}\otimes \mathbb{G}$,  i.e. 
  $$(\mathbb{F}\otimes \mathbb{G})_k=\mathop\oplus\limits_{i+j=k} \mathbb{F}_i\otimes \mathbb{G}_j.$$
    Then,
  $\Tor_1^S(S/I,S/J)=0$ implies that $\mathbb{F}\otimes \mathbb{G}$ is the minimal free resolution of $S/(I+J)$. 
  Let $$P:=\mathop\sum\limits_{i,j}\beta_{i,i+j}(S/I)x^iy^{i+j}$$
 and 
$$Q:=\mathop\sum\limits_{i,j}\beta_{i,i+j}(S/J)x^iy^{i+j}$$
 be the graded Poincar\'e series of $S/I$ and resp. $S/J.$ Then $PQ$ is the graded Poincar\'e series of $S/(I+J).$
  Since 
  $$\reg(S/I)=\deg_y(P),\ \reg(S/J)=\deg_y(Q)$$
  and
  $$\reg(S/(I+J))=\deg_y(PQ)=\deg_y(P)+\deg_y(Q),$$ 
  we get that 
  $$\reg(I+J)=\reg(S/(I+J))+1=\reg(S/I)+\reg(S/J)+1=$$
  $$=\reg(I)-1+\reg(J)-1+1=q+r-1.$$
  
\begin{Corollary}\label{cor_sum}
  $\reg (I_q+J_r)=q+r-1$.
\end{Corollary}

\begin{Lemma}
$\reg (I_qJ_r)=q+r$.  
\end{Lemma}

\textit{Proof}:  We have the exact sequence 
  \[
  0 \to I_q J_r\to I_q\oplus J_r\to I_q+J_r\to 0.
  \]
Therefore, by Corollary \ref{cor_sum} and  of \cite[Cor. 20.19]{Ei} we have 
\[
q+r\leq \reg(I_q J_r) \leq \max(\reg(I_q\oplus J_r),\ \reg (I_q+J_r)+1 ),  
\]
and the assertion follows.

\begin{Remark} 
  {\rm We observe that, when we consider the mixed product ideal 
  \[I_q J_r+ I_s J_t,\]
  we may fix $q<s$ and $t<r$, since otherwise one of the two summands contains the other.}   
\end{Remark}

\begin{Remark} {\rm Let $I$ be a monomial ideal and $G(I)=\{u_1,\ldots,u_q\}$ be the set of minimal generators of $I.$ Consider the exact sequence
\[0\to Z_1(I) \to R^q \stackrel{\phi}{\to} I\to 0,\]
where $Z_1(I)$ is the first syzygy of $I.$
We recall  that if we put 
 $$\sigma_{ij}=\frac{u_i}{(u_i,u_j)} e_{u_j} - \frac{u_i}{(u_i,u_j)} e_{u_i},\ 1\leq i,j\leq q, $$ 
 where $\phi(e_{u_i})=u_i$, $\phi(e_{u_j})=u_j$, then $\{\sigma_{ij}\}_{1\leq i,j\leq q}$ is a system  of generators of $Z_1(I)$ \cite[Lemma 15.1]{Ei}. For simplicity, we denote this system of generators of $Z_1(I)$ by $SG(I).$ }
\end{Remark}

\begin{Theorem}\label{theo_reg}
Suppose that $q<s$ and $t<r$. Then
\[\reg (I_qJ_r+I_sJ_t)=r+s-1.\]
\end{Theorem}

\textit{Proof}: We show first that $\reg (I_qJ_r+I_sJ_t)\leq r+s-1$.
In fact we have the exact sequence
  \[
  0 \to I_q J_r\cap I_sJ_t\to I_qJ_r\oplus I_sJ_t\to I_qJ_r+I_sJ_t\to 0.
  \]
  From the assumption we get
  $$ I_qJ_r\cap I_sJ_t=I_q\cap J_r\cap I_s \cap J_t=I_s J_r.$$
  Applying   \cite[Cor. 20.19]{Ei}  we obtain

 $$
    \reg (I_qJ_r+I_sJ_t)\leq \max (\reg (I_s J_r) -1, \reg (I_qJ_r\oplus I_sJ_t) )=$$
        $$                =   \max (s+r-1, q+r, s+t )=
                     r+s -1.$$
 To complete the proof it is sufficient to show that there exists an element $f\in Z_1(I_qJ_r+I_sJ_t)$,
  that has degree $r+s-1$ and is not generated by the other generators of $Z_1(I_qJ_r+I_sJ_t)$.

Since 
  $$G(I_qJ_r+I_sJ_t)=G(I_qJ_r)\cup G(I_sJ_t),$$
 we consider an element $f\in Z_1(I_qJ_r+I_sJ_t)$, 
  $$f=\frac{v}{(u,v)} e_u - \frac{u}{(u,v)} e_v,$$
  $$\ u=x_{i_1}\cdots x_{i_q}y_{j_1}\cdots y_{j_r}\in G(I_q J_r),$$
  $$ v=x_{i_1}\cdots x_{i_s}y_{j_1}\cdots y_{j_t}\in G(I_s J_t),$$
such that $i_1<i_2<\ldots<i_q<i_{q+1}<\ldots<i_s$ and $j_1<j_2<\ldots<j_t<j_{t+1}<\ldots<j_r$.
 Then
  $$f=x_{i_{q+1}}\cdots x_{i_s} e_u - y_{j_{t+1}}\cdots y_{j_r} e_v$$  
  and, if we consider the shift, the degree of $f$ is $r+s-1$. We observe also that $(u,v)$ has the minimal degree between the
  greatest common divisors of pairs of monomials respectively in $G(I_q J_r)$ and  $G(I_s J_t)$, and it cannot exist a syzygy of bigger degree that is a generator of $Z_1(I_qJ_r+I_sJ_t),$ since $\reg(I_qJ_r+I_sJ_t)\leq r+s-1.$ 
Let $g$ be a generator of $Z_1(I_qJ_r+I_sJ_t)$. Then either $g$ is in $SG(I_qJ_r)\cup
SG(I_rJ_s)$, or $g$ has degree  $r+s-1$ and is defined in a similar way as $f$.
If there exists $g$ of degree $r+s-1$  we are done. If there does not exist $g$ of degree $r+s-1,$ $f$ has to be generated by elements in $Z_1(I_q J_r)\subset Z_1(I_q J_r+I_s J_t)$ and elements in $Z_1(I_sJ_t)\subset Z_1(I_q J_r+I_s J_t)$. Since $Z_1(I_q J_r)\cap Z_1(I_s J_t)=\emptyset $ we obtain that 
  $$x_{i_{q+1}}\cdots x_{i_s} e_u\in Z_1(I_q J_r)
 \ {\rm and}\ y_{j_{t+1}}\cdots y_{j_r} e_v\in Z_1(I_s J_t),$$ that is false.

\section{ Depth and height of ideals of mixed products}
\label{sec:depth}

In this section we want to study the Krull dimension and the depth of the ring $S/(I_q J_r+I_s J_t)$. As a consequence we obtain conditions for the Cohen-Macaulayness of this ring. We start with an easy result.

\begin{Proposition}\label{Dimension}

\begin{enumerate}
\item If  $q=0,$ then $\dim(S/J_r)=n+r-1.$
\item If $r=0,$ then $\dim(S/I_q)=m+q-1.$
\item  If $q>0$ and $r>0,$ then 
 $$\dim (S/I_qJ_r)=n+m-\min(n-q+1,m-r+1).$$
 \end{enumerate}
\end{Proposition}

 \textit{Proof}: 1) and 2) follow immediately from Remark \ref{rem:mixdual}.
  \par\noindent 3) We have to show that 
  $$\height(I_q J_r)=\min(n-q+1,m-r+1).$$
   We observe that $I_q$, resp. $J_r$, is an unmixed
  ideal, whose associated prime ideals have height $n-q+1$, resp. $m-r+1$.
  The prime decomposition of the radical ideal $I_q J_r$ is 
  $$I_q J_r=I_q\cap J_r=(\mathop\bigcap\limits_{i=1}^t P_i)\cap(\mathop\bigcap\limits_{j=1}^s Q_j)$$
  where $\Ass(S/I_q)=\{P_1,\dots,P_t\}$ and $\Ass(S/J_r)=\{Q_1,\dots,Q_s\}$. The assertion follows.

\begin{Theorem} \label{dimensiune}
  Let $1\leq q<s$, $1\leq t<r$. Then
\begin{enumerate}
\item  $\dim S/(I_s+J_r)=r+s-2$;
  \item  $\dim S/(I_qJ_r+I_s)=n+m-\min(n-q+1,m+n-(r+s)+2)$;
  \item $\dim S/(I_qJ_r+I_s J_t)=n+m-\min(n-q+1,m-t+1,m+n-(r+s)+2)$.
\end{enumerate}
\end{Theorem}

\textit{Proof}: 1) We have  
$$S/(I_s+J_r)\cong K[x_1,\ldots,x_n]/I_s\otimes_K K[y_1,\ldots,y_m]/J_r$$ 
and the assertion follows by Remark \ref{rem:mixdual} and by \cite[Ex. 2.1.14]{Vi}.

2) Since $I_q \supset I_s,$ we have 
$$I_q J_r+I_s=I_q\cap (J_r+I_s).$$
 The assertion follows from the fact that each  face ideal of $J_r+I_s$ is the disjoint union of a face ideal of $I_s$ with a face ideal of $J_r$.

3) We remark that 
\[I_qJ_r+I_sJ_t=I_q\cap (J_r+I_s)\cap J_t\]
$I_q \supset I_s$, $J_t \supset J_r$. Now we continue as in case 2).



\begin{Lemma}\label{3.2}
  $\depth (S/I_qJ_r)\geq  q+r-1$.    
\end{Lemma}

\textit{Proof}: We consider the exact sequence
  
  $$0\to (I_q + J_r)/J_r \to S/J_r \to S/(I_q + J_r)\to 0.$$ 
  By  \cite[Prop. 1.2.9]{BH} we have 
  \[\depth (I_q + J_r)/J_r\geq \min( \depth (S/J_r), \depth( S/(I_q+J_r)+1)). \]
  By Remark \ref{rem:mixdual} and  Lemma \ref{lem_reg} we obtain that $\depth (S/J_r)=n+r-1$.
We also observe that 
  $$S/(I_q+J_r)\cong K[x_1,\ldots,x_n]/I_q  \otimes_K K[y_1,\ldots,y_m]/J_r.$$ 
 Then from \cite[Th. 2.2.21]{Vi} we get
  
    $$\depth S/(I_q+J_r)=\depth(K[x_1,\ldots,x_n]/I_q) +\depth(K[y_1,\ldots,y_m]/J_r)=$$
    $$=q-1+r-1=q+r-2.$$
  Therefore 
  $$\depth((I_q + J_r)/J_r)\geq \min(n+r-1,q+r-1)=q+r-1.$$
  If we consider the exact sequence 
  \[0\to I_q/I_q J_r \to S/I_q J_r \to S/I_q \to 0 \] 
 again by \cite[Prop. 1.2.9]{BH} we obtain 
  $$\depth(S/I_qJ_r) \geq \min(\depth(I_q/I_q J_r), \depth(S/I_q) ),$$
  and since 
  $$I_q/I_qJ_r=I_q/I_q\cap J_r\cong (I_q + J_r)/J_r,$$ we get 
  $$\depth(S/I_qJ_r) \geq \min(q+r-1,m+q-1)=q+r-1.$$

\begin{Theorem}\label{Depth}
\begin{enumerate}
\item If $q=0$, then $\depth(S/J_r)=n+r-1.$
\item If $r=0,$ then $\depth(S/I_q)=m+q-1.$
\item Suppose that $q>0$ and $r>0.$ Then  
$$\depth (S/I_qJ_r)=q+r-1.$$
\end{enumerate}
\end{Theorem}

\textit{Proof}: 1) and 2) follow immediately from the proof of \ref{lem_reg}.
\par\noindent 3) It remains to show that  $\depth(S/I_qJ_r)\leq q+r-1$. We prove this by induction on $r$ and $q$. 
Let $r=1$, $q=1.$ We claim that 
  \[\Tor_{m+n-1}(K,S/I_1J_1)\neq 0\]
  which implies the inequality.
  We shall construct an element  
  \[[z]\in H:=H_{m+n-1}(K(S;\mathbf{x},\mathbf{y})\otimes S/I_1J_1),\]
  where $K(S;\mathbf{x},\mathbf{y})$ is the Koszul complex of $S$ with respect to the sequence $x_1,\ldots,x_n,y_1,\ldots,y_m$.

  We fix the lexicographic term ordering such that
  $$x_1>\cdots>x_n>y_1>\cdots>y_m$$
   and let 
  \[z:=\sum_{k=1}^m (-1)^{k+1}y_k (\bigwedge_{i=1}^n e_i)\wedge (\bigwedge_{j\neq k} f_j).\]
 We want to show that $z$ is a cycle.
 $$\partial(z)=\sum_{k=1}^m (-1)^{k+1} y_k\partial[(\bigwedge_{i=1}^n  e_i)\wedge (\bigwedge_{j\neq k} f_j)]$$
    $$=\sum_{k=1}^m (-1)^{k+1} y_k  [ (-1)^n (\bigwedge_{i=1}^n e_i) \partial(\bigwedge_{j\neq k} f_j)] +  
    \sum_{k=1}^m (-1)^{k+1}y_k[ \partial( \bigwedge_{i=1}^n e_i )\bigwedge_{j\neq k} f_j].$$ 
It is easy to see that the second sum vanishes, since  $x_iy_k \in I_1J_1$, for all $i=1,\dots,n,\ k=1,\dots,m.$ Therefore
  \[\partial(z)=(-1)^n(\sum_{k=1}^m (-1)^{k+1}y_k  \bigwedge_{i=1}^n e_i \partial(\bigwedge_{j\neq k}
  f_j)),\] hence we have to show that
  \begin{equation}
    \label{eq:cycle}
    \sum_{k=1}^m (-1)^{k+1} y_k  \partial(\bigwedge_{j\neq k} f_j)=0.
  \end{equation}
 A summand of  (\ref{eq:cycle}) is of the form
    \[(-1)^{k+1} y_k (-1)^{\pi(l)+1} y_{l}(\bigwedge_{j\neq k,\ l}f_j),\] 
    where $\pi(l)$ is the position of $f_l$ in $(\bigwedge_{j\neq k} f_j)$. We observe that if $k<l$ then $\pi(l)=l-1$, while if
    $k>l$ then $\pi(l)=l$.
 We also observe that  in (\ref{eq:cycle}) there exists only one other summand 
    \[(-1)^{k'+ \pi(l')} y_{k'} y_{l'}(\bigwedge_{j\neq k',\ l'}f_j), \]
    with $l'=k$ and $k'=l$. It is easy to show that 
    \[(-1)^{l+ \pi(k)} y_{l} y_{k} + (-1)^{k+\pi(l)}y_k y_{l}=0\]
    Indeed, either $k<l$, that is $\pi(k)=k$ and  $\pi(l)=l-1$, or $k>l$, that is  $\pi(k)=k-1$ and $\pi(l)=l$.

    Now we want to show that $z$ is not a boundary. We observe that, if we write 
    \[z= (e_1\wedge \ldots\wedge e_n)\wedge( \sum_{k=1}^m (-1)^{k+1} y_k (\bigwedge_{j\neq k} f_j)),\]
    \begin{enumerate}
    \item The coefficients of $z$ are polynomials of total degree 1; 
    \item The multidegree of the coefficients of $z$ with respect to $\mathbf{x}$ is 
      \[\multideg_{\mathbf{x}}(z)=(0,\ldots,0)\in {\bf Z}^n.\]
    \end{enumerate}
    Suppose now that $b$ is a boundary such that $\partial (b)=z$. Then
   $$b=\gamma e_1\wedge \ldots \wedge e_n \wedge f_1 \wedge \ldots\wedge f_m,\ \gamma\in S$$
    and
   $$\partial(b)= \gamma \sum_{k=1}^n (-1)^{k+1} x_k (\bigwedge_{i\neq k} e_i)\wedge (\mathop\bigwedge\limits_{j=1}^m f_j)+$$ 
 $$+ \gamma \sum_{l=1}^m (-1)^{n+l+1} y_l (\mathop\bigwedge\limits_{i=1}^n e_i)\wedge (\bigwedge_{j\neq l} f_j).$$

   Suppose that  $\gamma \neq 0$. From the fact that the multidegree of $z$ with respect to $\mathbf{x}$ is $(0,\ldots,0)$, we get that 
   $\gamma x_i= 0$ for all $i=1,\dots,n$. It follows that $\gamma x_i\in I_1 J_1$ and the total degree of  $\gamma$ is at least 1. 
Therefore $\gamma y_l$ has total degree at least 2, for all $l=1,\dots,m$. But the coefficients of $z$ have total degree $1,$ therefore $\gamma=0$.

    Now suppose that $1<i\leq m$ and that $\depth(S/I_1 J_{i-1})=i-1$. We want to show that $\depth(S/I_1 J_{i})=i.$

    Let $S_l=K[x_1,\ldots,x_n,y_1,\ldots, y_l]$, with $l=1,\ldots,m$ and let 
    $$L_i:=I_1 J_i\subseteq S_l.$$
    $$L^+_j:=I_1 J_j\subseteq S_{l+1}.$$

    We have the exact sequence 
    \[0\to L_i^+/L_iS_{l+1}\to S_{l+1}/L_iS_{l+1}\to S_{l+1}/L^+_i\to 0.\]
    We observe that 
    $$L_i^+=L_iS_{l+1}+y_{l+1}L_{i-1}S_{l+1}$$
     and
    $$L_i^+/L_iS_{l+1}= (L_iS_{l+1}+y_{l+1}L_{i-1}S_{l+1})/L_i\cong$$
    $$\cong y_{l+1} L_{i-1}S_{l+1}/L_iS_{l+1}\cap y_{l+1}L_{i-1}S_{l+1}.$$
    Since $L_i\subset L_{i-1},\ y_{l+1}\cap L_{i-1}S_{l+1}=y_{l+1}L_{i-1}S_{l+1}$ and $y_{l+1}$ is not a zero-divisor on $S_{l+1},$ we obtain
    $$L_i^+/L_iS_{l+1}\cong L_{i-1}S_{l+1}/L_iS_{l+1}.$$ 

    Therefore 
    $$S_{l+1}/L_i^+\cong S_{l+1}/L_{i-1}S_{l+1}.$$ 
    By the induction hypothesis we have
    $$\depth (S_l/L_{i-1})=i-1$$
     and since $y_{l+1}$
    is a regular element of $S_{l+1}/L_{i-1}S_{l+1}$ we obtain the assertion.
If we consider the ideal $I_iJ_r$, with the same argument, by induction on $i\geq 1$, we obtain that $\depth (S/I_qJ_r)=q+r-1$.

\begin{Corollary}\label{Coh}
$S/(I_qJ_r)$ is Cohen-Macaulay if and only if one of the following conditions hold:
\begin{enumerate}
\item $q=0;$
\item $r=0;$
\item $q>0, r>0, m=r$ and $n=q.$
\end{enumerate}
\end{Corollary}

\textit{Proof}: The cases $q=0$ or $r=0$ being clear, we can suppose that $q>0$ and $ r>0.$ By \ref{Dimension} we have 
$$\dim (S/I_qJ_r)=m+n-\min(n-q+1,m-r+1)$$
 and by \ref{Depth} we get
$$\depth (S/I_qJ_r)=q+r-1.$$
It follows that if $S/(I_qJ_r)$ is Cohen-Macaulay then 
$$m+n-\min(n-q+1,m-r+1)=q+r-1.$$
Suppose that $n-q+1\leq m-r+1.$ Then from above we get  $m=r$ and then clearly $n=q$. The case $m-r+1\leq n-q+1$ leads to the same condition. The converse is obvious. 

\begin{Lemma}\label{le:Depth}
  Let $1\leq q<s$ and $1\leq t<r$. Then
\begin{enumerate}
 \item $\depth S/(I_qJ_r+I_sJ_t)\geq  \min (q+r-1, s+t-1 );$
 \item $\depth S/(I_qJ_r+I_s)\geq q+r-1.$
\end{enumerate}
\end{Lemma}

\textit{Proof}:  1) We have the exact sequence 
  \[0\to S/I_qJ_r\cap I_sJ_t \to S/I_qJ_r\oplus S/I_sJ_t\to S/(I_qJ_r+I_sJ_t)\to 0\]
  where $I_qJ_r\cap I_sJ_t=I_sJ_r$ (see the proof of Theorem \ref{theo_reg}).
  By \cite[Prop. 1.2.9]{BH}  we obtain 
  $$\depth (S/(I_qJ_r+I_sJ_t))\geq \min ( \depth(S/I_sJ_r)-1, \depth((S/I_qJ_r)\oplus(S/I_sJ_t)) )=$$
   $$= \min ( s+r-2,\min(q+r-1, s+t-1 )) =$$
   $$= \min (q+r-1, s+t-1 ).$$

2) In this case we have
\[0\to S/I_sJ_r \to S/I_qJ_r\oplus S/I_s\to S/(I_qJ_r+I_s)\to 0\]
and we obtain
$$\depth (S/(I_qJ_r+I_s))\geq \min ( \depth(S/I_sJ_r)-1, \depth((S/I_qJ_r)\oplus(S/I_s)) )=$$
   $$= \min ( s+r-2,\min(q+r-1, s+m-1 )) = q+r-1.$$

\begin{Theorem} \label{3333}
  Let $1\leq q<s$, $1\leq t<r$. Then:
  \begin{enumerate}
  \item $\depth S/(I_s+J_r)=s+r-2$;
  \item $\depth S/(I_qJ_r+I_s)=q+r-1;$ 
  \item $\depth S/(I_qJ_r+I_s J_t)=\min(q+r,s+t)-1.$
\end{enumerate}
\end{Theorem}
\textit{Proof}: 1) This was already proved in the proof of Lemma \ref{3.2}.
\par\noindent  2) Since by Lemma \ref{le:Depth} we have that $\depth S/(I_qJ_r+I_s )\geq q+r-1,$ we only have to prove that
\[\Tor_{m+n-(q+r-1)}(K,S/(I_qJ_r+I_s))\neq 0.\]
 We consider the exact sequence
  \begin{equation}
    \label{eq:exact1}
    0\to S/I_s J_r\to S/I_qJ_r\oplus S/I_s\to S/(I_qJ_r+I_s )\to 0.
  \end{equation}
  
  Let $k=n+m-(q+r-1)$ and consider the long exact sequence of $\Tor$ induced by the sequence (\ref{eq:exact1}). We get the exact sequence
  \begin{eqnarray*}
    \cdots\to\Tor_k(K,S/I_s J_r)&\to&\Tor_k(K,S/I_qJ_r\oplus S/I_s)\stackrel{\phi}{\to}\\
    &\stackrel{\phi}{\to}& \Tor_k(K,S/(I_qJ_r+I_s)) \to\cdots
  \end{eqnarray*}
   By hypothesis we have $s>q$ and by Theorem \ref{Depth} we get 
   $$\pd S/I_sJ_r=m+n-(s+r-1)<k.$$
    Therefore we have
  that $\Tor_k(K,S/I_s J_r)=0$, that is $\phi$ is injective.
 We have also that 
  \[\Tor_k(K,S/I_qJ_r\oplus S/I_sJ_t)\cong \Tor_k(K,S/I_qJ_r)\oplus  \Tor_k(K,S/I_sJ_t)\] 
  and by Theorem \ref{Depth} we obtain $\Tor_k(K,S/I_qJ_r)\neq 0$.
 Therefore by the injectivity of $\phi$ we obtain the assertion.

\par\noindent  3) Let us suppose   that $q+r\leq s+t$. Since by Lemma \ref{le:Depth} we have that $\depth S/(I_qJ_r+I_s J_t)\geq q+r-1$
  we have to show that $\depth S/(I_qJ_r+I_s J_t)\leq q+r-1$, that is equivalent to say
  \[\Tor_{m+n-(q+r-1)}(K,S/(I_qJ_r+I_s J_t))\neq 0.\]
  We consider the exact sequence
  \begin{equation}
    \label{eq:exact}
    0\to S/I_s J_r\to S/I_qJ_r\oplus S/I_sJ_t\to S/(I_qJ_r+I_s J_t)\to 0.
  \end{equation}
  Let $k=m+n-(q+r-1)$ and consider the long exact sequence of $\Tor$ induced by sequence (\ref{eq:exact}). We get
  \begin{eqnarray*}
    \cdots\to\Tor_k(K,S/I_s J_r)&\to&\Tor_k(K,S/I_qJ_r\oplus S/I_sJ_t)\stackrel{\psi}{\to}\\
    &\stackrel{\psi}{\to}& \Tor_k(K,S/(I_qJ_r+I_sJ_t)) \to\cdots
  \end{eqnarray*}
  By hypothesis we have $s>q$ and by Theorem \ref{Depth} we get 
  $$\pd S/I_sJ_r=m+n-(s+r-1)<k.$$
Therefore we have
  that $\Tor_k(K,S/I_s J_r)=0$, that is $\psi$ is injective.
We have also that 
  \[\Tor_k(K,S/I_qJ_r\oplus S/I_sJ_t)\cong \Tor_k(K,S/I_qJ_r)\oplus  \Tor_k(K,S/I_sJ_t)\] 
  and by Theorem \ref{Depth} we obtain $\Tor_k(K,S/I_qJ_r)\neq 0$.
 Therefore by the injectivity of $\psi$ we obtain the assertion.

\begin{Corollary}\label{CM} Let $0\leq q<s$, $0\leq t<r$. Then:
\begin{enumerate}
\item $S/I_q, S/J_r$ and $S/(I_q+J_r)$ are always Cohen-Macaulay;
\item $S/(I_qJ_r)$ is Cohen-Maculay if and only if $m=r$ and $n=q$;
\item $S/(I_qJ_r+I_s)$ is Cohen-Macaulay if and only if $s=q+1$ and $r=m$;
\item $S/(I_qJ_r+I_sJ_t)$ is Cohen-Macaulay if and only if $r=m, s=n,\\ t=m-1,  q=n-1$.
\end{enumerate}
\end{Corollary}

\textit{Proof}: (1) and (2) see Corollary \ref{Coh}. 

\noindent (3) Suppose that $R:=S/(I_qJ_r+I_s)$ is Cohen-Macaulay. Then from \ref{3333} and \ref{dimensiune} we get 
$$n+m-\min(n-q+1, n+m-(r+s)+2)=q+r-1.$$ 
We have 2 cases:

a) $m+q\geq r+s-1$. Then $\dim(R)=m+q-1$ and $\depth(R)=q+r-1.$ Since $R$ is Cohen-Macaulay we get $m=r$. From this it follows that $s\leq q+1$ and consequently $s=q+1.$

b) $m+q<r+s-1$. Then $\dim(R)=r+s-2$. Since $R$ is Cohen-Macaulay it follows at once that $s=q+1$. Further we get $m<r$, which is a contradiction.
\par\noindent The converse is obvious.

\par\noindent (4) Let $R:=S/(I_qJ_r+I_sJ_t).$ and suppose that $R$ is Cohen-Macaulay. Then from \ref{3333} and \ref{dimensiune} we have that 
$$\min(q+r,s+t)-1=n+m-\min(n-q+1,m-t+1,n+m-(r+s)+2).$$
We have  to consider 2 cases:
\par a) $q+r\leq s+t.$ We shall consider 3 subcases, as follows.
\par a1) $n-q+1\leq m-t+1$ and $n-q+1\leq n+m-r-s+2.$ 
\par Since $R$ is Cohen-Macaulay we obtain that $m=r.$  From 
$$n-q\leq n+m-r-s+1$$ it follows that $q\geq s-1$, that is $q=s-1.$ Now from $$s-1+r\leq s+t$$
 we get $r-1\leq t$ and consequently $$t=r-1=m-1.$$ Further we obtain immediately $n=q+1.$
\par a2) $n+m-r-s+2\leq m-t+1$ and $n+m-r-s+2\leq n-q+1.$ 
This means that 
$$m+q\leq r+s-1,$$ $$n+t\leq r+s-1.$$
Since $R$ is Cohen-Macaulay we get immediately that $q=s-1$ and from the second inequality above we have $r\geq m$, that is $r=m.$ From $q+r\leq s+t$ we have that 
$t\geq r-1$, that is $t=m-1.$ Now it follows at once that $s=n.$

\par a3) $m-t+1\leq n-q+1 $ and $m-t+1\leq n+m-r-s+2.$ 
\par Since $R$ is Cohen-Macaulay we get $q+r=n+t.$ Then 
$$r=n-q+t\geq m-t+t=m,$$
 that is $r=m.$ From 
$$m-t\leq m+n-r-s+1$$ we get 
$$q+r=t+n\geq r+s-1.$$
It follows at once that $q\geq s-1$, that is $q=s-1.$  Now, from 
$$q+r=n+t\leq s+t$$ we obtain $s=n$ and then $t=q+r-n=s-1+r-s=r-1.$
\par b) $s+t\leq q+r.$ One has to consider the same 3 subcases as in case a). The proof is exactly the same.
\par\noindent The converse is clear.

 \vspace{1.0cm}


\end{document}